# PROPERTIES OF PRINCIPAL COMPONENT METHODS FOR FUNCTIONAL AND LONGITUDINAL DATA ANALYSIS

By Peter Hall, Hans-Georg Müller [1] and Jane-Ling Wang[2]

*Australian National University, University of California, Davis and University of California, Davis*

The use of principal component methods to analyze functional data is appropriate in a wide range of different settings. In studies of "functional data analysis," it has often been assumed that a sample of random functions is observed precisely, in the continuum and without noise. While this has been the traditional setting for functional data analysis, in the context of longitudinal data analysis a random function typically represents a patient, or subject, who is observed at only a small number of randomly distributed points, with nonnegligible measurement error. Nevertheless, essentially the same methods can be used in both these cases, as well as in the vast number of settings that lie between them. How is performance affected by the sampling plan? In this paper we answer that question. We show that if there is a sample of $n$ functions, or subjects, then estimation of eigenvalues is a semiparametric problem, with root-$n$ consistent estimators, even if only a few observations are made of each function, and if each observation is encumbered by noise. However, estimation of eigenfunctions becomes a nonparametric problem when observations are sparse. The optimal convergence rates in this case are those which pertain to more familiar function-estimation settings. We also describe the effects of sampling at regularly spaced points, as opposed to random points. In particular, it is shown that there are often advantages in sampling randomly. However, even in the case of noisy data there is a threshold sampling rate (depending on the number of functions treated) above which the rate of sampling (either randomly or regularly) has negligible impact on estimator performance, no matter whether eigenfunctions or eigenvectors are being estimated.

Received August 2004; revised July 2005.

[1]Supported in part by NSF Grants DMS-03-54448 and DMS-05-05537.

[2]Supported in part by NSF Grant DMS-04-06430.

## 1. Introduction.

1.1. *Connections between FDA and LDA.* Advances in modern technology, including computing environments, have facilitated the collection and analysis of high-dimensional data, or data that are repeated measurements of the same subject. If the repeated measurements are taken over a period of time, say on an interval $\mathscr{I}$, there are generally two different approaches to treating them, depending on whether the measurements are available on a dense grid of time points, or whether they are recorded relatively sparsely.

When the data are recorded densely over time, often by machine, they are typically termed functional or curve data, with one observed curve (or function) per subject. This is often the case even when the data are observed with experimental error, since the operation of smoothing data recorded at closely spaced time points can greatly reduce the effects of noise. In such cases we may regard the entire curve for the $i$th subject, represented by the graph of the function $X_i(t)$ say, as being observed in the continuum, even though in reality the recording times are discrete. The statistical analysis of a sample of $n$ such graphs is commonly termed functional data analysis, or FDA, and can be explored as suggested in the monographs by Ramsay and Silverman [27, 28].

Biomedical longitudinal studies are similar to FDA in important respects, except that it is rare to observe the entire curve. Measurements are often taken only at a few scattered time points, which vary among subjects. If we represent the observation times for subject $i$ by random variables $T_{ij}$, for $j = 1, \ldots, m_i$, then the resulting data are $(X_i(T_{i1}), \ldots, X_i(T_{im_i}))$, generally observed with noise. The study of information in this form is often referred to as longitudinal data analysis, or LDA. See, for example, [12] or [20].

Despite the intrinsic similarities between sampling plans for functional and longitudinal data, statistical approaches to analyzing them are generally distinct. Parametric technologies, such as generalized estimating equations or generalized linear mixed effects models, have been the dominant methods for longitudinal data, while nonparametric approaches are typically employed for functional data. These and related issues are discussed by Rice [31].

A significant, intrinsic difference between the two settings lies in the perception that functional data are observed in the continuum, without noise, whereas longitudinal data are observed at sparsely distributed time points and are often subject to experimental error. However, functional data are sometimes computed after smoothing noisy observations that are made at a relatively small number of time points, perhaps only a dozen points if, for example, full-year data curves are calculated from monthly figures (see, e.g., [26].) Such instances indicate that the differences between the two data types relate to the way in which a problem is perceived and are arguably more conceptual than actual—for example, in the case of FDA, as one where discretely recorded data are more readily understood as observations of a continuous process.

As this discussion suggests, in view of these close connections, there is a need to understand the interface between FDA and LDA views of data that might reasonably be thought of as having a functional origin. This is one of the goals of the present paper. In the context of principal component analysis, we explain the effect that observation at discrete time points, rather than observation in the continuum, has on statistical estimators. In particular, and as we shall show, estimators of the eigenvalues $\theta_j$ of principal components can be root-$n$ consistent even when as few as two observations are made of each of the $n$ subjects, and even if experimental error is present. However, in such cases, estimation of eigenfunctions $\psi_j$ is at slower rates, but nevertheless at rates that would be optimal for function estimators if data on those functions were observed in conventional form. On the other hand, when the $n$ random functions are fully observed in the continuum, the convergence rates of both eigenvalue and eigenfunction estimators are $n^{-1/2}$.

These results can be summarized by stating that estimation of $\theta_j$ or of $\psi_j$ is a semiparametric problem when the random functions are fully observed in the continuum, but that estimation of $\psi_j$ is a nonparametric problem, whereas estimation of $\theta_j$ remains semiparametric, when data are observed sparsely with noise. Indeed, if the number of observations per subject is at least two but is bounded, and if the covariance function of subjects has $r$ bounded derivatives, then the minimax-optimal, mean square convergence rate of eigenfunction estimators is $n^{-2r/(2r+1)}$. This rate is achieved by estimators based on empirical spectral decomposition. We shall treat in detail only the case $r = 2$, since that setting corresponds to estimation of covariance using popular local smoothing methods. However, straightforward arguments give the extension to general $r$.

We also identify and discuss the important differences between sampling at regularly spaced, and at random time points. Additionally we address the case where the number of sampled points per subject increases with sample size. Here we show that there is a threshold value rate of increase which ensures that estimators of $\theta_j$ and of $\psi_j$ are first-order equivalent to their counterparts in the case where subjects are fully observed, without noise. By drawing connections between FDA, where nonparametric methods are well developed and popular, and LDA, where parametric techniques play a dominant role, we are able to point to ways in which nonparametric methodology may be introduced to LDA. There is a range of settings in LDA where parametric models are difficult to postulate. This is especially true when the longitudinal data do not have similar "shapes," or are so sparse that the individual data profiles cannot be discerned. Measurement errors can also mask the shapes of the underlying subject profiles. Thus, more flexible models based on nonparametric approaches are called for, at least at the early stage of data analysis. This motivates the application of functional data approaches, and in particular, functional principal component analysis, to longitudinal data.

It might be thought that our analysis of the infinite-dimensional problem of FDA should reveal the same phenomena that are apparent in "large $p$, small $n$" theory for finite-dimensional problems. For example, estimators of the maximum eigenvalue might be expected to be asymptotically biased. See, for example, Johnstone [18]. However, these features are not present in theoretical studies of conventional FDA methodology (see, e.g., [4, 11]), where complete functions are observed, and it is arguably unsurprising that they are not present. One reason is that, although FDA is infinite-dimensional, an essential ingredient distinguishing it from the multivariate vector case is smoothness. The problem of estimating any number of eigenvalues and eigenfunctions does not become successively more difficult as sample size increases, since this problem in some sense may be reduced to that of estimating fixed smooth mean and covariance functions from the available functional data. In contrast, in typical "large $p$, small $n$" asymptotics, the dimension of covariance matrices is assumed to increase with sample size which gives rise to specific properties.

The results in the present paper represent the first attempt at developing concise asymptotic theory and optimal rates of convergence describing functional PCA for sparse data. Upper bounds for rates of convergence of estimated eigenfunctions in the sparse-data case, but not attaining the concise convergence rates given in the present paper, were developed by Yao, Müller and Wang [38] under more restrictive assumptions. Other available theoretical results for functional PCA are for the ideal situation when entire random functions are observed, including Dauxois, Pousse and Romain [11], Bosq [3], Pezzulli and Silverman [25], Boente and Fraiman [2], Cardot [7], Girard [14] and Hall and Hosseini-Nasab [15]. There is an extensive literature on the general statistical analysis of functional data when the full functions are assumed known. It includes work of Besse and Ramsay [1], Castro, Lawton and Sylvestre [10], Rice and Silverman [32], Brumback and Rice [5] and Cardot, Ferraty and Sarda [8, 9], as well as many articles discussed and cited by Ramsay and Silverman [27, 28]. Kneip and Utikal [21] used methods of functional data analysis to assess the variability of densities for data sets from different populations. Contributions to various aspects of the analysis of sparse functional data, including longitudinal data observed with measurement error, include those of Shi, Weiss and Taylor [34], Staniswalis and Lee [35], James, Hastie and Sugar [17], Rice and Wu [33] and Müller [24]. For practical issues of implementing and applying functional PCA, we refer to [6, 19, 29, 32, 37, 38].

## 2. Functional PCA for discretely observed random functions.

2.1. *Functional principal component analysis.* Let $X_1, \ldots, X_n$ denote independent and identically distributed random functions on a compact interval $\mathscr{I}$, satisfying $\int_{\mathscr{I}} E(X^2) < \infty$. The mean function is $\mu = E(X)$, and the covariance function is $\psi(u, v) = \operatorname{cov}\{X(u), X(v)\}$. Functional PCA is based on interpreting $\psi$ as the kernel of a linear mapping on the space $L_2(\mathscr{I})$ of square-integrable functions on $\mathscr{I}$, taking $\alpha$ to $\psi\alpha$ defined by $(\psi\alpha)(u) = \int_{\mathscr{I}} \alpha(v)\psi(u, v)\, dv$. For economy we use the same notation for an operator and its kernel. Mercer's theorem (e.g., [16], Chapter 4) now implies a spectral decomposition of the function $\psi$,

$$(2.1) \qquad \psi(u, v) = \sum_{j=1}^{\infty} \theta_j \psi_j(u)\psi_j(v),$$

where $\theta_1 \geq \theta_2 \geq \cdots \geq 0$ are ordered values of the eigenvalues of the operator $\psi$, and the $\psi_j$'s are the corresponding eigenfunctions.

The eigenfunctions form a complete orthonormal sequence on $L_2(\mathscr{I})$, and so we may represent each function $X_i - \mu$ in terms of its generalized Fourier expansion in the $\psi_j$'s,

$$(2.2) \qquad X_i(t) = \mu(t) + \sum_{j=1}^{\infty} \zeta_{ij}\psi_j(t),$$

where $\zeta_{ij} = \int_{\mathscr{I}} (X_{i-\mu})\psi_j$ is referred to as the $j$th functional principal component score, or random effect, of the $i$th subject, whose (observed as in FDA or hidden as in LDA) random trajectory is $X_i$. The expansion (2.2) is referred to as the Karhunen–Loève or functional principal component expansion of the stochastic process $X_i$. The fact that $\psi_j$ and $\psi_k$ are orthogonal for $j \neq k$ implies that the random variables $\zeta_{ij}$, for $1 \leq j < \infty$, are uncorrelated.

Although the convergence in (2.2) is in $L_2$, not pointwise in $t$, the only purpose of that result, from the viewpoint of this paper, is to define the principal components, or individual effects, $\zeta_{ij}$. The values of those random variables are defined by (2.2), with probability 1.

The difficulty of representing distributions of random functions means that principal component analysis assumes even greater importance in the setting of FDA than it does in more conventional, finite-dimensional statistical problems. Especially if $j$ is small, the shape of the function $\psi_j$ conveys interpretable information about the shapes one would be likely to find among the curves in the data set $X_1, \ldots, X_n$, if the curves were observable. In particular, if $\psi_1$ has a pronounced turning point in a part of $\mathscr{I}$ where the other functions, with relatively low orders, are mostly flat, then the turning point is likely to appear with high probability in a random function $X_i$. The "strength" with which this, or another, feature is likely to arise is proportional to the standard deviation of $\zeta_{ij}$, that is, to the value of $\theta_j^{1/2}$. Conversely, if all eigenfunctions are close to zero in a given region, we may conclude that the underlying random process is constrained to be close to its mean in this region with relatively little random variation. These considerations and others, including the fact that (2.1) can be used to represent a variety of characteristics of the random function $X$, motivate the development of methods for estimating each $\theta_j$ and each $\psi_j$.

2.2. *Estimation.* Let $X_1, \ldots, X_n$ be as in Section 1.1, and assume that for each $i$ we observe pairs $(T_{ij}, Y_{ij})$ for $1 \leq j \leq m_i$, where

$$(2.3) \qquad Y_{ij} = X_i(T_{ij}) + \varepsilon_{ij},$$

the "observation times" $T_{ij}$ all lie in the compact interval $\mathscr{I}$, the errors $\varepsilon_{ij}$ have zero mean and each $m_j \geq 2$. For simplicity, when developing theoretical properties it will be assumed that the $T_{ij}$'s are identically distributed random variables, that the errors $\varepsilon_{ij}$ are also identically distributed, with finite variance $E(\varepsilon^2) = \sigma^2$, and that the $X_i$'s, $T_{ij}$'s and $\varepsilon_{ij}$'s are totally independent. However, similar results may be obtained with a degree of weak dependence, and in cases of nonidentical distributions.

Using the data set $\mathcal{D} = \{(T_{ij}, Y_{ij}), \ 1 \leq j \leq m_i, \ 1 \leq i \leq n\}$, we wish to construct estimators $\hat{\theta}_j$ and $\hat{\psi}_j$ of $\theta_j$ and $\psi_j$, respectively. We start with estimators $\hat{\mu}$ of $\mu = E(X)$, and $\hat{\psi}$ of the autocovariance, $\psi$; definitions of $\hat{\mu}$ and $\hat{\psi}$ will be given shortly. The function $\hat{\psi}$, being symmetric, enjoys an empirical version of the expansion at (2.1),

$$(2.4) \qquad \hat{\psi}(u, v) = \sum_{j=1}^{\infty} \hat{\theta}_j \hat{\psi}_j(u) \hat{\psi}_j(v).$$

Here, $\hat{\theta}_1, \hat{\theta}_2, \ldots$ are eigenvalues of the operator $\hat{\psi}$, given by $(\hat{\psi}\alpha)(u) = \int_{\mathscr{I}} \alpha(v) \times \hat{\psi}(u, v) \, dv$ for $\alpha \in L_2(\mathscr{I})$, and $\hat{\psi}_j$ is the eigenfunction corresponding to $\hat{\theta}_j$. In Section 3 we shall develop properties of $\hat{\theta}_j$ and $\hat{\psi}_j$. Given $j_0 \geq 1$, the $\hat{\theta}_j$'s are ordered so that $\hat{\theta}_1 \geq \cdots \geq \hat{\theta}_{j_0} \geq \hat{\theta}_j$, the last inequality holding for all $j > j_0$.

The signs of $\psi_j$ and $\hat{\psi}_j$ can be switched without altering either (2.1) or (2.4). This does not cause any difficulty, except that, when discussing the closeness of $\psi_j$ and $\hat{\psi}_j$, we clearly want them to have the same parity. That is, we would like these eigenvectors to "point in the same general direction" when they are close. We ensure this by allowing the sign of $\psi_j$ to be chosen arbitrarily, but asking that the sign of $\hat{\psi}_j$ be chosen to minimize $\|\hat{\psi}_j - \psi_j\|$ over both possible choices, where here and in the following $\|\cdot\|$ denotes the $L_2$-norm, $\|\psi\| = (\int \psi^2)^{1/2}$.

We construct first $\widehat{\mu}(u)$, and then $\widehat{\psi}(u,v)$, by least-squares fitting of a local linear model, as follows. Given $u \in \mathscr{I}$, let $h_\mu$ and $h_\phi$ denote bandwidths and select $(\hat{a}, \hat{b}) = (a, b)$ to minimize

$$\sum_{i=1}^{n} \sum_{j=1}^{m_i} \{Y_{ij} - a - b(u - T_{ij})\}^2 K\left(\frac{T_{ij} - u}{h_\mu}\right),$$

and take $\widehat{\mu}(u) = \hat{a}$. Then, given $u, v \in \mathscr{I}$, choose $(\hat{a}_0, \hat{b}_1, \hat{b}_2) = (a_0, b_1, b_2)$ to minimize

$$\sum_{i=1}^{n} \sum_{j,k\,:\,1 \leq j \neq k \leq m_i} \{Y_{ij}Y_{ik} - a_0 - b_1(u - T_{ij}) - b_2(v - T_{ik})\}^2$$

$$\times K\left(\frac{T_{ij} - u}{h_\phi}\right) K\left(\frac{T_{ik} - v}{h_\phi}\right).$$

The quantity $\hat{a}_0$ estimates $\phi(u,v) = E\{X(u)X(v)\}$, and so we denote it by $\widehat{\phi}(u,v)$. Put

$$\widehat{\psi}(u,v) = \widehat{\phi}(u,v) - \widehat{\mu}(u)\widehat{\mu}(v).$$

These estimates are the same as those proposed in [38], where practical features regarding the implementation are discussed in detail. The emphasis in [38] is on estimating the random effects $\zeta_{ij}$, for which Gaussian assumptions are made on processes and errors. We extend the consistency results for eigenvalues and eigenvectors of [38] in four significant ways: First, we establish concise first-order properties. Second, the first-order results in the present paper imply bounds that are an order of magnitude smaller than the upper bounds provided by Yao, Müller and Wang [38]. Third, we derive the asymptotic distribution of estimated eigenvalues. Fourth, we characterize a transition where the asymptotics of "longitudinal behavior" with sparse and noisy measurements per subject transform into those of "functional behavior" where random trajectories are completely observed. This transition occurs as the number of measurements per subject is allowed to increase at a certain rate.

The operator defined by $\widehat{\psi}$ is not, in general, positive semidefinite, and so the eigenvalues $\hat{\theta}_j$ at (2.4) may not all be negative. Nevertheless, $\widehat{\psi}$ is symmetric, and so (2.4) is assured.

Define $U_{ij} = u - T_{ij}$, $V_{ik} = v - T_{ik}$, $Z_{ijk} = Y_{ij}Y_{ik}$,

$$W_{ij} = K\left(\frac{T_{ij} - u}{h_\mu}\right), \qquad W_{ijk} = K\left(\frac{T_{ij} - u}{h_\phi}\right) K\left(\frac{T_{ik} - v}{h_\phi}\right).$$

Using this notation we may write

(2.5) $\quad \widehat{\mu}(u) = \dfrac{S_2 R_0 - S_1 R_1}{S_0 S_2 - S_1^2}, \qquad \widehat{\phi}(u,v) = (A_1 R_{00} - A_2 R_{10} - A_3 R_{01})B^{-1},$

where

$$A_1 = S_{20}S_{02} - S_{11}^2, \qquad A_2 = S_{10}S_{02} - S_{01}S_{11}, \qquad A_3 = S_{01}S_{20} - S_{10}S_{11},$$

$$S_r = \sum_{i=1}^{n}\sum_{j=1}^{m_i} U_{ij}^r W_{ij}, \qquad\qquad R_r = \sum_{i=1}^{n}\sum_{j=1}^{m_i} U_{ij}^r Y_{ij} W_{ij},$$

$$S_{rs} = \sum_{i,j,k:\,j<k}\sum\sum U_{ij}^r V_{ik}^s W_{ijk}, \qquad R_{rs} = \sum_{i,j,k:\,j<k}\sum\sum U_{ij}^r V_{ik}^s Z_{ijk} W_{ijk},$$

$$B = A_1 S_{00} - A_2 S_{10} - A_3 S_{01}$$

$$= \left\{\sum_{i,j,k:\,j<k}\sum\sum (U_{ij} - \bar U)^2 W_{ijk}\right\}\left\{\sum_{i,j,k:\,j<k}\sum\sum (V_{ij} - \bar V)^2 W_{ijk}\right\}$$

$$- \left\{\sum_{i,j,k:\,j<k}\sum\sum (U_{ij} - \bar U)(V_{ij} - \bar V) W_{ijk}\right\}^2 \geq 0,$$

$$\bar Q = \left(\sum_{i,j,k:\,j<k}\sum\sum Q_{ij} W_{ijk}\right)\Big/\left(\sum_{i,j,k:\,j<k}\sum\sum W_{ijk}\right),$$

for $Q = U, V$. Here we have suppressed the dependence of $S_r$, $R_r$ and $W_{ij}$ on $u$, and of $A_r$, $B$, $S_{rs}$, $R_{rs}$ and $W_{ijk}$ on $(u, v)$.

## 3. Theoretical properties.

3.1. *Main theorems.* Our estimators $\hat\mu$ and $\hat\psi$ have been constructed by local linear smoothing, and so it is natural to make second derivative assumptions below, as a prerequisite to stating both upper and lower bounds to convergence rates. If $\hat\mu$ and $\hat\psi$ were defined by $r$th-degree local polynomial smoothing, then we would instead assume $r$ derivatives, and in particular the optimal $L_2$ convergence rate of $\hat\psi_j$ would be $n^{-r/(2r+1)}$ rather than the rate $n^{-2/5}$ discussed below.

Assume that the random functions $X_i$ are independent and identically distributed as $X$ and are independent of the errors $\varepsilon_{ij}$; that the latter are independent and identically distributed as $\varepsilon$, with $E(\varepsilon) = 0$ and $E(\varepsilon^2) = \sigma^2$; that

$$\text{for each } C > 0 \qquad \max_{j=0,1,2} E\left\{\sup_{u \in \mathscr{I}} |X^{(j)}(u)|^C\right\} + E(|\varepsilon|^C) < \infty;$$

that the kernel function $K$ is compactly supported, symmetric and Hölder continuous; that for an integer $j_0 > 1$ there are no ties among the $j_0 + 1$ largest eigenvalues of $\phi$ [although we allow the $(j_0 + 1)$st largest eigenvalue to be tied with the $(j_0 + 2)$nd]; that the data pairs $(T_{ij}, Y_{ij})$ are observed for $1 \leq j \leq m_i$ and $1 \leq i \leq n$, where each $m_i \geq 2$ and $\max_{i \leq n} m_i$ is bounded as $n \to \infty$; that the $T_{ij}$'s have a common distribution, the density, $f$, of which is bounded away from zero on $\mathcal{I}$; and that $n^{\eta - (1/2)} \le h_\mu = o(1)$, for some $\eta > 0$. In addition, for parts (a) and (b) of Theorem 1 we assume, respectively, that (a) $n^{\eta - (1/3)} \le h_\phi$ for some $\eta > 0$, $\max(n^{-1/3} h_\phi^{2/3}, n^{-1} h_\phi^{-8/3}) = o(h_\mu)$, $h_\mu = o(h_\phi)$ and $h_\phi = o(1)$; and (b) $n^{\eta - (3/8)} \le h_\phi$ and $h_\phi + h_\mu = o(n^{-1/4})$. Call these conditions (C).

In conditions (C) above we suppose the $m_i$'s to be deterministic, but with minor modifications they can be taken to be random variables. Should some of the $m_i$'s be equal to 1, these values may be used to estimate the mean function, $\mu$, even though they cannot contribute to estimates of the covariance. For simplicity we shall not address this case, however.

Put $x = X - \mu$, and define $\kappa = \int K^2$, $\kappa_2 = \int u^2 K(u)\, du$,

$$c(r, s) = \int_{\mathcal{I}} f(t)^{-1} \psi_r(t) \psi_s(t)\, dt,$$

$$(3.1) \quad \beta(u, v, w, z) = E\{x(u)x(v)x(w)x(z)\} - \psi(u, v)\psi(w, z),$$

$$\chi(u, v) = \tfrac{1}{2} \kappa_2 \{\psi_{20}(u, v) + \psi_{02}(u, v) + \mu''(u)\mu(v) + \mu(u)\mu''(v)\},$$

where $\psi_{rs}(u, v) = (\partial^{r+s}/\partial u^r\, \partial v^s)\psi(u, v)$. Let

$$N = \tfrac{1}{2} \sum_{i \le n} m_i(m_i - 1)$$

and

$$\nu(r, s) = \sum_{i=1}^{n} \sum_{j_1 < k_1} \sum_{j_2 < k_2} E[\{f(T_{ij_1})f(T_{ik_1})f(T_{ij_2})f(T_{ik_2})\}^{-1}$$
$$\times \beta(T_{ij_1}, T_{ik_1}, T_{ij_2}, T_{ik_2})$$
$$\times \psi_r(T_{ij_1})\psi_r(T_{ik_1})\psi_s(T_{ij_2})\psi_s(T_{ik_2})].$$

This formula has a conceptually simpler, although longer to write, version, obtained by noting that the $T_{ij}$'s are independent with density $f$. Asymptotic bias and variance properties of estimators are determined by the quantities

$$C_1 = C_1(j) = \kappa \iint_{\mathcal{I}^2} \frac{E\{x(t_1)^2 x(t_2)^2\} + \sigma^2}{f(t_1)f(t_2)} \psi_j(t_1)^2\, dt_1\, dt_2,$$

$$(3.2) \quad C_2 = C_2(j) = \sum_{k:k \ne j} (\theta_j - \theta_k)^{-2} \left(\int \chi \psi_j \psi_k\right)^2,$$

$$(\Sigma)_{rs} = N^{-2}\{\nu(r, s) + N\sigma^2 c(r, s)^2\}.$$

Let $\Sigma$ denote the $j_0 \times j_0$ matrix with $(r, s)$th component equal to $(\Sigma)_{rs}$. Note that $(\Sigma)_{rs} = O(n^{-1})$ for each pair $(r, s)$. Write $\vec{a}$, $\vec{\theta}$ and $\vec{\hat{\theta}}$ for the vectors $(a_1, \ldots, a_{j_0})^{\mathrm{T}}$, $(\theta_1, \ldots, \theta_{j_0})^{\mathrm{T}}$ and $(\hat{\theta}_1, \ldots, \hat{\theta}_{j_0})^{\mathrm{T}}$, respectively.

Our next result describes large-sample properties of eigenvalue and eigenfunction estimators. It is proved in Section 4.

THEOREM 1.   *Assume conditions* (C). *Then,* (a) *for* $1 \leq j \leq j_0$,

$$(3.3) \qquad \|\widehat{\psi}_j - \psi_j\|^2 = \frac{C_1}{Nh_\phi} + C_2 h_\phi^4 + o_p\{(nh_\phi)^{-1} + h_\phi^4\},$$

*and* (b) *for any vector* $\vec{a}$, $\vec{a}^{\mathrm{T}}(\widehat{\vec{\theta}} - \vec{\theta})$ *is asymptotically normally distributed with mean zero and variance* $\vec{a}^{\mathrm{T}} \Sigma \vec{a}$.

The representation in part (b) of the theorem is borrowed from [13]. Bounds on $\widehat{\psi}_j - \psi_j$ and on $\widehat{\theta}_j - \theta_j$, which hold uniformly in increasingly large numbers of indices $j$, and in particular for $1 \leq j \leq j_0 = j_0(n)$ where $j_0(n)$ diverges with $n$, can also be derived. Results of this nature, where the whole functions $X_i$ are observed without error, rather than noisy observations being made at scattered times $T_{ij}$ as in (2.3), are given by Hall and Hosseini-Nasab [15]. The methods there can be extended to the present setting. However, in practice there is arguably not a great deal of interest in moderate- to large-indexed eigenfunctions. As Ramsay and Silverman [28] note, it can be difficult to interpret all but relatively low-order principal component functions.

The order of magnitude of the right-hand side of (3.3) is minimized by taking $h \asymp n^{-1/5}$; the relation $a_n \asymp b_n$, for positive numbers $a_n$ and $b_n$, means that $a_n/b_n$ is bounded away from zero and infinity as $n \to \infty$. Moreover, it may be proved that if $h \asymp n^{-1/5}$, then the relation $\|\widehat{\psi}_j - \psi_j\| = O_p(n^{-2/5})$, implied by (3.3), holds uniformly over a class of distributions of processes $X$ that satisfy a version of conditions (C). The main interest, of course, lies in establishing the reverse inequality, uniformly over all candidates $\widetilde{\psi}_j$ for estimators of $\psi_j$, thereby showing that the convergence rate achieved by $\widehat{\psi}_j$ is minimax-optimal.

We shall do this in the case where only the first $r$ eigenvalues $\theta_1, \ldots, \theta_r$ are nonzero, with fixed values, where the joint distribution of the Karhunen–Loève coefficients $\zeta_{i1}, \ldots, \zeta_{ir}$ [see (2.2)] is also fixed, and where the observation times $T_{ij}$ are uniformly distributed. These restrictions actually strengthen the lower bound result, since they ensure that the "max" part of the minimax bound is taken over a relatively small number of options.

The class of eigenfunctions $\psi_j$ will be taken to be reasonably rich, however; we shall focus next on that aspect. Given $c_1 > 0$, let $\mathcal{S}(c_1)$ denote the $L_\infty$ Sobolev space of functions $\phi$ on $\mathscr{I}$ for which $\max_{s=0,1,2} \sup_{t \in \mathscr{I}} |\phi^{(s)}(t)| \leq c_1$. We pass from this space to a class $\Psi = \Psi(c_1)$ of vectors $\vec{\psi} = (\psi_1, \ldots, \psi_r)$ of orthonormal functions, as follows. Let $\psi_{11}, \ldots, \psi_{1r}$ denote any fixed functions that are orthonormal on $\mathscr{I}$ and have two bounded, continuous derivatives there. For each sequence $\phi_1, \ldots, \phi_r$ of functions in $\mathcal{S}(c_1)$, let $\psi_1, \ldots, \psi_r$ be the functions constructed by applying Gram–Schmidt orthonormalization to $\psi_{11} + \phi_1, \ldots, \psi_{1r} + \phi_r$, working through this sequence in any order but nevertheless taking $\psi_j$ to be the new function obtained on adjoining $\psi_{1j} + \phi_j$ to the orthonormal sequence, for $1 \le j \le r$. If $c_1$ is sufficiently small, then, for some $c_2 > 0$,

$$\sup_{\vec{\psi} \in \Psi} \max_{1 \le j \le r} \max_{s=0,1,2} \sup_{t \in \mathscr{I}} |\psi_j^{(s)}(t)| \le c_2.$$

Moreover, defining $\mathcal{A}_j$ to be the class of functions $\psi_{2j} = \psi_{1j} + \phi_j$ for which $\phi_j \in \mathcal{S}(c_1)$ and $\int \psi_{2j}^2 = 1$, we have

$$(3.4) \qquad \mathcal{A}_j(c_1) \subseteq \{\psi_j : (\psi_1, \ldots, \psi_r) \in \Psi\}$$

for $1 \le j \le r$. In the discussion below we shall assume that these properties hold.

Let $\theta_1 > \cdots > \theta_r > 0$ be fixed, and take $0 = \theta_{r+1} = \theta_{r+2} = \cdots$. Let $\zeta_1, \ldots, \zeta_r$ be independent random variables with continuous distributions, all moments finite, zero means, and $E(\zeta_j^2) = \theta_j$ for $1 \le j \le r$. Assume we observe data $Y_{ij} = X_i(T_{ij}) + \varepsilon_{ij}$, where $1 \le i \le n$, $1 \le j \le m$, $m \ge 2$ is fixed, $X_i = \sum_{1 \le k \le r} \zeta_{ik} \psi_j$, $(\psi_1, \ldots, \psi_r) \in \Psi$, each $\vec{\zeta_i} = (\zeta_{i1}, \ldots, \zeta_{ir})$ is distributed as $(\zeta_1, \ldots, \zeta_r)$, each $T_{ij}$ is uniformly distributed on $\mathscr{I} = [0, 1]$, the $\varepsilon_{ij}$'s are identically normally distributed with zero mean and nonzero variance, and the $\vec{\zeta_i}$'s, $T_{ij}$'s and $\varepsilon_{ij}$'s are totally independent. Let $\widetilde{\Psi}_j$ denote the class of all measurable functionals $\widetilde{\psi}_j$ of the data $\mathcal{D} = \{(T_{ij}, Y_{ij}), 1 \le i \le n, 1 \le j \le m\}$. Theorem 2 below asserts the minimax optimality in this setting of the $L_2$ convergence rate $n^{-2/5}$ for $\widehat{\psi}_j$ given by Theorem 1.

THEOREM 2. *For the above prescription of the data* $\mathcal{D}$, *and assuming* $h \asymp n^{-1/5}$,

$$(3.5) \qquad \lim_{C \to \infty} \limsup_{n \to \infty} \max_{1 \le j \le r} \sup_{\vec{\psi} \in \Psi} P(\|\widehat{\psi}_j - \psi_j\| > C n^{-2/5}) = 0;$$

*and for some* $C > 0$,

$$(3.6) \qquad \liminf_{n \to \infty} \min_{1 \le j \le r} \inf_{\widetilde{\psi}_j \in \widetilde{\Psi}_j} \sup_{\vec{\psi} \in \Psi} P\{\|\widetilde{\psi}_j(\mathcal{D}) - \psi_j\| > C n^{-2/5}\} > 0.$$

It is possible to formulate a version of (3.5) where, although the maximum over $j$ continues to be in the finite range $1 \le j \le r$, the supremum over $\vec{\psi} \in \Psi$ is replaced by a supremum over a class of infinite-dimensional models. There one fixes $\theta_1 > \cdots > \theta_r > \theta_{r+1} \ge \theta_{r+2} \ge \cdots \ge 0$, and chooses $\psi_{r+1}, \psi_{r+2}, \ldots$ by extension of the process used to select $\psi_1, \ldots, \psi_r$.

3.2. *Discussion.* Part (b) of Theorem 1 gives the asymptotic joint distribution of the components $\hat{\theta}_j$, and from part (a) we may deduce that the asymptotically optimal choice of $h_\phi$ for estimating $\psi_j$ is $h_\phi \sim (C_1/4C_2N)^{1/5} \asymp n^{-1/5}$. More generally, if $h_\phi \asymp n^{-1/5}$, then by Theorem 1 $\|\hat{\psi}_j - \psi_j\| = O_p(n^{-2/5})$. By Theorem 2 this convergence rate is asymptotically optimal under the assumption that $\psi$ has two derivatives. For $h_\phi$ of size $n^{-1/5}$, the conditions on $h_\mu$ imposed for part (a) of Theorem 1 reduce to $n^{-7/15} = o(h_\mu)$ and $h_\mu = o(n^{-1/5})$.

In particular, Theorem 1 argues that a degree of undersmoothing is necessary for estimating $\psi_j$ and $\theta_j$. Even when estimating $\psi_j$, the choice of $h_\phi$ can be viewed as undersmoothed, since the value const. $n^{-1/5}$, suggested by (3.3), is an order of magnitude smaller than the value that would be optimal for estimating $\phi$; there the appropriate size of $h_\phi$ is $n^{-1/6}$. The suggested choice of $h_\mu$ is also an undersmoothing choice, for estimating both $\psi_j$ and $\theta_j$.

Undersmoothing, in connection with nonparametric nuisance components, is known to be necessary in situations where a parametric component of a semiparametric model is to be estimated relatively accurately. Examples include the partial-spline model studied by Rice [30], and extensions to longitudinal data discussed by Lin and Carroll [22]. In our functional PCA problem, where the eigenvalues and eigenfunctions are the primary targets, the mean and covariance functions are nuisance components. The fact that they should be undersmoothed reflects the cases mentioned just above, although the fact that one of the targets is semiparametric, and the other nonparametric, is a point of departure.

The assumptions made about the $m_i$'s and $T_{ij}$'s in Theorems 1 and 2 are realistic for sparsely sampled subjects, such as those encountered in longitudinal data analysis. There, the time points $T_{ij}$ typically represent the dates of biomedical follow-up studies, where only a few follow-up visits are scheduled for each patient, and at time points that are convenient for that person. The result is a small total number of measurements per subject, made at irregularly spaced points.

On the other hand, for machine-recorded functional data the $m_i$'s are usually larger and the observation times are often regularly spaced. Neither of the two theorems is valid if the observation times $T_{ij}$ are of this type, rather than (as in the theorems) located at random points. For example, if each $m_i = m$ and we observe each $X_i$ only at the points $T_{ij} = j/m$, for $1 \leq j \leq m$, then we cannot consistently estimate either $\theta_j$ or $\psi_j$ from the resulting data, even if no noise is present. There exist infinitely many distinct distributions of random functions $X$, in particular of Gaussian processes, for which the joint distribution of $X(j/m)$, for $1 \leq j \leq m$, is common. The stochastic nature of the $T_{ij}$'s, which allows them to take values arbitrarily close to any given point in $\mathscr{I}$, is critical to Theorems 1 and 2. Nevertheless, if the $m_i$'s increase sufficiently quickly with $n$, then regular spacing is not a problem. We shall discuss this point in detail shortly.

Next we consider how the results reported in Theorem 1 alter when each $m_i \geq m$ and $m$ increases. However, we continue to take the $T_{ij}$'s to be random variables, considered to be uniform on $\mathcal{I}$ for simplicity. In this setting, $N \geq \frac{1}{2}m(m-1)n$ and $\nu(r,s) = \frac{1}{4}m^4 nd(r,s) + O(m^3n)$, where

$$d(r,s) = \int \beta(u,v,w,z)\psi_r(u)\psi_r(v)\psi_s(w)\psi_s(z)\,du\,dv\,dw\,dz.$$

If we assume each $m_i = m$, then it follows that $(\Sigma)_{rs} = n^{-1}d(r,s) + O\{(mn)^{-1}\}$, as $m, n \to \infty$. The leading term here, that is, $n^{-1}d(r,s)$, equals the $(r,s)$th component of the limiting covariance matrix of the conventional estimator of $(\theta_1,\ldots,\theta_{j_0})^{\mathrm{T}}$ when the full curves $X_i$ are observed without noise. [It may be shown that the proof leading to part (b) of Theorem 1 remains valid in this setting, where each $m_i = m$ and $m = m(n) \to \infty$ as $n$ increases.] This reflects the fact that the noisy sparse-data, or LDA, estimators of eigenvalues converge to their no-noise and full-function, or FDA, counterparts as the number of observations per subject increases, no matter how slow the rate of increase.

The effect on $\widehat{\psi}_j$ of increasing $m$ is not quite as clear from part (a) of Theorem 1. That result implies only that $\|\widehat{\psi}_j - \psi_j\|^2 = C_2 h_\phi^4 + o_p\{(nh_\phi)^{-1}\}$. Therefore the order of magnitude of the variance component is reduced, and a faster $L_2$ convergence rate of $\widehat{\psi}_j$ to $\psi_j$ can be achieved by choosing $h_\phi$ somewhat smaller than before. However, additional detail is absent.

To obtain further information it is instructive to consider specifically the "FDA approach" when full curves are not observed. There, a smooth function estimator $\widehat{X}_i$ of $X_i$ would be constructed by passing a statistical smoother through the sparse data set $\mathcal{D}_i = \{(T_{ij}, Y_{ij}), 1 \leq j \leq m_i\}$, with $Y_{ij}$ given by (2.3). Functional data analysis would then proceed as though $\widehat{X}_i$ were the true curve, observed in its entirety. The step of constructing $\widehat{X}_i$ is of course a function estimation one, and should take account of the likely smoothness of $X_i$. For example, if each $X_i$ had $r$ derivatives, then a local polynomial smoother of degree $r-1$ might be passed through $\mathcal{D}_i$. (Kneip and Utikal [21] also employed a conventional smoother, this time derived from kernel density estimation, in their exploration of the use of functional-data methods for assessing different population densities.)

Let us take $r = 2$ for definiteness, in keeping with the assumptions leading to Theorems 1 and 2, and construct $\widehat{X}_i$ by running a local-linear smoother through $\mathcal{D}_i$. Assume that each $m_i \geq m$. Then the following informally stated property may be proved: The smoothed function estimators $\widehat{X}_i$ are as good as the true functions $X_i$, in the sense that the resulting estimators of both $\theta_j$ and $\psi_j$ are first-order equivalent to the root-$n$ consistent estimators that arise on applying conventional principal component analysis to the true curves $X_i$, provided $m$ is of larger order than $n^{1/4}$. Moreover, this result holds true for both randomly distributed observation times $T_{ij}$ and regularly spaced times. A formal statement and outline proof of this result are given in Section 3.4.

These results clarify issues that are sometimes raised in FDA and LDA, about whether effects of "the curse of dimensionality" have an impact through the number of observations per subject. It can be seen from our results that having $m_i$ large is a blessing rather than a curse; even in the presence of noise, statistical smoothing successfully exploits the high-dimensional character of the data and fills in the gaps between adjacent observation times.

Theorem 1 provides advice on how the bandwidth $h_\phi$ might be varied for different eigenfunctions $\psi_j$. It suggests that, while the order of magnitude of $h_\phi$ need not depend on $j$, the constant multiplier could, in many instances, be increased with $j$. The latter suggestion is indicated by the fact that, while the constant $C_1$ in (3.3) will generally not increase quickly with $j$, $C_2$ will often tend to increase relatively quickly, owing to the spacings between neighboring eigenvalues decreasing with increasing $j$. The connection to spacings is mathematically clear from (3.2), where it is seen that by decreasing the values of $\theta_j - \theta_k$ we increase $C_2$. Operationally, it is observed that higher-order empirical eigenfunctions are typically increasingly oscillatory, and hence require more smoothing for effective estimation.

3.3. *Proof of Theorem* 2. The upper bound (3.5) may be derived using the argument in Section 4. To obtain (3.6) it is sufficient to show that if $j \in [1, r]$ is fixed and the orthonormal sequence $\{\psi_1, \ldots, \psi_r\}$ is constructed starting from $\psi_j \in \mathcal{A}_j(c_1)$; and if, in addition to the data $\mathcal{D}$, the values of each $\zeta_{ik}$, for $1 \le i \le n$ and $1 \le k \le m$, and of each $\psi_k(T_{i\ell})$, for $1 \le i \le n$, $k \ne j$ and $1 \le \ell \le m$, are known; then for some $C_j > 0$,

$$\liminf_{n \to \infty} \inf_{\tilde{\psi}_j \in \tilde{\Psi}_j} \sup_{\psi_j \in \mathcal{A}_j(c_1)} P\{\|\tilde{\psi}_j(\mathcal{D}) - \psi_j\| > C_j n^{-2/5}\} > 0.$$

[To obtain this equivalence we have used (3.4).] That is, if we are given only the data $\mathcal{D}' = \{(T_{ij}, \zeta_{ij}, \psi_j(T_{ij}) + \varepsilon_{ij}\zeta_{ij}^{-1}),\ 1 \le i \le n,\ 1 \le j \le m\}$, and if $\bar{\Psi}_j$ denotes the class of measurable functions $\bar{\psi}_j$ of $\mathcal{D}'$, then it suffices to show that for some $C_j > 0$,

$$\liminf_{n \to \infty} \inf_{\bar{\psi}_j \in \bar{\Psi}_j} \sup_{\psi_j \in \mathcal{A}_j(c_1)} P\{\|\bar{\psi}_j(\mathcal{D}') - \psi_j\| > C_j n^{-2/5}\} > 0.$$

Except for the fact that the errors here are $\varepsilon_{ij}\zeta_{ij}^{-1}$ rather than simply $\varepsilon_{ij}$, this result is standard; see, for example, [36]. The factor $\zeta_{ij}^{-1}$ is readily dealt with by using a subsidiary argument.

3.4. *Random function approximation.* In Section 3.2 we discussed an approach to functional PCA that was based on running a local-linear smoother through increasingly dense, but noisy, data on the true function $X_i$, producing an empirical approximation $\widehat{X}_i$. Here we give a formal statement and outline proof of the result discussed there.

THEOREM 3. *Suppose each $m_i \geq m$, and assume conditions* (C) *from Section 3.1, except that the observation times $T_{ij}$ might be regularly spaced on $\mathcal{I}$ rather than being randomly distributed there. Estimate $\theta_j$ and $\psi_j$ using conventional PCA for functional data, as though each smoothed function estimator $\widehat{X}_i$ really were the function $X_i$. Then the resulting estimators of $\theta_j$ and $\psi_j$ are root-n consistent, and first-order equivalent to the conventional estimators that we would construct if the $X_i$'s were directly observed, provided $m = m(n)$ diverges with $n$ and the bandwidth, $h$, used for the local-linear smoother satisfies $h = o(n^{-1/4})$, $mhn^{-\delta_1} \to \infty$ and $m^{1-\delta_2}h \to \infty$, for some $\delta_1, \delta_2 > 0$.*

We close with a proof. Observe that the estimator of $\psi(u, v)$ that results from operating as though each $\widehat{X}_i$ is the true function $X_i$, is

$$\check{\psi}(u, v) = \frac{1}{n} \sum_{i=1}^{n} \{\widehat{X}_i(u) - \overline{\widehat{X}}(u)\}\{\widehat{X}_i(v) - \overline{\widehat{X}}(v)\},$$

where $\overline{\widehat{X}} = n^{-1} \sum_i \widehat{X}_i$. The linear operator that is defined in terms of $\check{\psi}$ is positive semidefinite. The FDA-type estimators $\check{\theta}_j$ and $\check{\psi}_j$ of $\theta_j$ and $\psi_j$, respectively, would be constructed by simply identifying terms in the corresponding spectral expansion,

$$\check{\psi}(u, v) = \sum_{j=1}^{\infty} \check{\theta}_j \check{\psi}_j(u) \check{\psi}_j(v),$$

where $\check{\theta}_1 \geq \check{\theta}_2 \geq \cdots \geq 0$. Of course, if we were able to observe the process $X_i$ directly, without noise, we would estimate $\psi$ using

$$\bar{\psi}(u, v) = \frac{1}{n} \sum_{i=1}^{n} \{X_i(u) - \bar{X}(u)\}\{X_i(v) - \bar{X}(v)\},$$

where $\bar{X} = n^{-1} \sum_i X_i$, and take as our estimators of $\theta_j$ and $\psi_j$ the corresponding terms $\bar{\theta}_j$ and $\bar{\psi}_j$ in the expansion,

$$\bar{\psi}(u, v) = \sum_{j=1}^{\infty} \bar{\theta}_j \bar{\psi}_j(u) \bar{\psi}_j(v),$$

with $\bar{\theta}_1 \geq \bar{\theta}_2 \geq \cdots \geq 0$.

Methods used to derive limit theory for the estimators $\bar{\theta}_j$ and $\bar{\psi}_j$ (see, e.g., [15]) may be used to show that the estimator pairs $(\hat{\theta}_j, \hat{\psi}_j)$ and $(\bar{\theta}_j, \bar{\psi}_j)$ are asymptotically equivalent to first order if $\check{\psi} - \hat{\psi} = o_p(n^{-1/2})$, but generally not first-order equivalent if $\check{\psi}$ and $\hat{\psi}$ differ in terms of size $n^{-1/2}$ or larger. Here, distances are measured in terms of the conventional $L_2$ metric for functions. Since we have used a local-linear smoother to construct the functions $\hat{X}_i$ from the data $\mathcal{D}$, then the bias contribution to $\check{\psi} - \hat{\psi}$ is of size $h^2$, where $h$ denotes the bandwidth for the local-linear method. The contribution from the error about the mean is of size $(mnh)^{-1/2}$ at each fixed point. The "penalty" to be paid for extending uniformly to all points is smaller than any polynomial in $n$. Indeed, using an approximation on a lattice that is of polynomial fineness, the order of magnitude of the uniform error about the mean can be seen to be of order $n^\delta(mnh)^{-1/2}$ for each $\delta > 0$. Theorem 3 follows.

## 4. Proof of Theorem 1.

*Step* (i): *Approximation lemmas.* Let $\psi$ denote a symmetric, strictly positive-definite linear operator on the class $L_2(\mathscr{I})$ of square-integrable functions from the compact interval $\mathscr{I}$ to the real line, with a kernel, also denoted by $\psi$, having spectral decomposition given by (2.1). Denote by $\bar{\psi}$ another symmetric, linear operator on $L_2(\mathscr{I})$. Write the spectral decomposition of $\bar{\psi}$ as $\bar{\psi}(u,v) = \sum_{j\geq 1} \bar{\theta}_j \bar{\psi}_j(u)\bar{\psi}_j(v)$. Since $\psi$ is nonsingular, then its eigenfunctions $\psi_j$, appearing at (2.1), comprise a complete orthonormal sequence, and so we may write $\bar{\psi}_j = \sum_{k\geq 1} \bar{a}_{jk}\psi_k$, for constants $\bar{a}_{jk}$ satisfying $\sum_{k\geq 1} \bar{a}_{jk}^2 = 1$. We may choose $\bar{a}_{jj}$ to be either positive or negative, since altering the sign of an eigenfunction does not change the spectral decomposition. Below, we adopt the convention that each $\bar{a}_{jj} \geq 0$.

Given a function $\alpha$ on $\mathscr{I}^2$, define $\|\alpha\| = (\iint_{\mathscr{I}^2} \alpha^2)^{1/2}$, $\|\alpha\|_\infty = \sup |\alpha|$ and

$$\|\alpha\|_{(j)}^2 = \int_{\mathscr{I}} \left\{ \int_{\mathscr{I}} \alpha(u,v)\psi_j(v)\,dv \right\}^2 du.$$

If $\alpha_1$ and $\alpha_2$ are functions on $\mathscr{I}$, write $\int \alpha \alpha_1 \alpha_2$ to denote

$$\iint_{\mathscr{I}^2} \alpha(u,v)\alpha_1(u)\alpha_2(v)\,du\,dv.$$

For example, $\int (\bar{\psi} - \psi)\psi_j\psi_j$ in (4.2), below, is to be interpreted in this way. Let $\int \alpha \alpha_1$ denote the function of which the value at $u$ is $\int_{\mathscr{I}} \alpha(u,v)\alpha_1(v)\,dv$, and write $|\mathscr{I}|$ for the length of $\mathscr{I}$.

LEMMA 1.  *For each* $j \geq 1$,

(4.1) $$\|\bar{\psi}_j - \psi_j\|^2 = 2(1 - \bar{a}_{jj}),$$

$$(4.2) \quad \begin{aligned} &\left| \bar{\theta}_j - \theta_j - \left\{ \int (\bar{\psi} - \psi) \psi_j \psi_j + (1 - \bar{a}_{jj}^2) \left( 1 - \int \bar{\psi} \psi_j \psi_j \right) \right\} \right| \\ &\quad \leq |\mathscr{I}| \|\bar{\psi}_j - \psi_j\|^2 \|\bar{\psi} - \psi\|_\infty \\ &\qquad + (1 - \bar{a}_{jj}^2) \|\bar{\psi} - \psi\| + 2\|\bar{\psi}_j - \psi_j\| \|\bar{\psi} - \psi\|_{(j)}. \end{aligned}$$

Lemma 1 implies that knowing bounds for $1 - \bar{a}_{jj}$ and for several norms of $\bar{\psi} - \psi$ gives us information about the sizes of $\|\bar{\psi}_j - \psi_j\|$ and $\bar{\theta}_j - \theta_j$. We shall take $\bar{\psi} = \hat{\psi}$, in which case we have an explicit formula for $\|\bar{\psi} - \psi\|$. Therefore our immediate need is for an approximation to $\bar{a}_{jj}$, denoted below by $\hat{a}_{jj}$ when $\bar{\psi} = \hat{\psi}$. This requirement will be filled by the next lemma. Define $\Delta = \hat{\psi} - \psi$, and let $\hat{a}_{jk}$ denote the generalized Fourier coefficients for expressing $\hat{\psi}_j$ in terms of the $\psi_k$'s: $\hat{\psi}_j = \sum_{k \geq 1} \hat{a}_{jk} \psi_k$, where we take $\hat{a}_{jj} \geq 0$.

LEMMA 2. *Under the conditions of Theorem 1,*

$$(4.3) \quad 1 - \hat{a}_{jj}^2 = O_p(\|\Delta\|_{(j)}^2),$$

$$(4.4) \quad \left| \hat{a}_{jj}^2 - 1 + \sum_{k : k \neq j} (\theta_j - \theta_k)^{-2} \left( \int \Delta \psi_j \psi_k \right)^2 \right| = O_p(\|\Delta\| \|\Delta\|_{(j)}^2).$$

Lemma 1 is derived by using basic manipulations in operator theory. The proof of Lemma 2 involves more tedious arguments, which can be considered to be sparse-data versions of methods employed by Bosq [4] to derive his Theorem 4.7 and Corollaries 4.7 and 4.8; see also [23].

*Step* (ii): *Implications of approximation lemmas.* Since $2(1 - \hat{a}_{jj}) = 1 - \hat{a}_{jj}^2 + O_p(|1 - \hat{a}_{jj}^2|^2)$ and $\|\Delta\|_{(j)} \leq \|\Delta\|$, then Lemma 2 implies that

$$(4.5) \quad 2(1 - \hat{a}_{jj}) = D_{j1} + O_p(\|\Delta\| \|\Delta\|_{(j)}^2),$$

where

$$D_{j1} = \sum_{k : k \neq j} (\theta_j - \theta_k)^{-2} \left( \int \Delta \psi_j \psi_k \right)^2 \leq \text{const.} \|\Delta\|_{(j)}^2.$$

Note too that

$$(4.6) \quad D_{j1} = D_{j2} + \theta_j^{-2} \|\Delta\|_{(j)}^2 - \theta_j^{-2} \left( \int \Delta \psi_j \psi_j \right)^2,$$

where

$$(4.7) \quad D_{j2} = \sum_{k : k \neq j} \left\{ (\theta_j - \theta_k)^{-2} - \theta_j^{-2} \right\} \left( \int \Delta \psi_j \psi_k \right)^2.$$

Standard arguments on uniform convergence of nonparametric function estimators can be used to show that, under the conditions of Theorem 1, $\|\widehat{\mu} - \mu\|_\infty = o_p(1)$ and $\|\widehat{\phi} - \phi\|_\infty = o_p(1)$, from which it follows that $\|\widehat{\psi} - \psi\|_\infty = o_p(1)$. Combining (4.5) and (4.6) with (4.1)–(4.4) we deduce that

$$
\begin{aligned}
(4.8) \qquad \|\widehat{\psi}_j - \psi_j\|^2 = {}& D_{j2} + \theta_j^{-2}\|\Delta\|^2_{(j)} - \theta_j^{-2}\Big(\int \Delta\psi_j\psi_j\Big)^2 \\
& + O_p(\|\Delta\|\,\|\Delta\|^2_{(j)}),
\end{aligned}
$$

$$
(4.9) \qquad \widehat{\theta}_j - \theta_j = \int \Delta\psi_j\psi_j + O_p(\|\Delta\|^2_{(j)}).
$$

Let $E'$ denote expectation conditional on the observation times $T_{ij}$, for $1 \le j \le m_i$ and $1 \le i \le n$. Standard methods may be used to prove that, under the bandwidth conditions imposed in either part of Theorem 1, and for each $\eta > 0$,

$$
E'\|\Delta\|^2 = O_p\{(nh_\phi^2)^{-1} + (nh_\mu)^{-1} + h^4\},
$$

and $E'\|\Delta\|^2_{(j)} = O_p\{(nh)^{-1} + h^4\}$, where $h = h_\phi + h_\mu$. Therefore, under the bandwidth assumptions made, respectively, for parts (a) and (b) of Theorem 1, the "$O_p$" remainder term on the right-hand side of (4.8) equals $o_p\{(nh_\phi)^{-1} + h_\phi^4\} + O_p\{(nh_\mu)^{-3/2} + h_\mu^6\}$, while the remainder on the right-hand side of (4.9) equals $o_p(n^{-1/2})$. Hence,

$$
\begin{aligned}
(4.10) \qquad \|\widehat{\psi}_j - \psi_j\|^2 = {}& D_{j2} + \theta_j^{-2}\|\Delta\|^2_{(j)} - \theta_j^{-2}\Big(\int \Delta\psi_j\psi_j\Big)^2 \\
& + o_p\{(nh_\phi)^{-1} + h_\phi^4\} + O_p\{(nh_\mu)^{-3/2} + h_\mu^6\},
\end{aligned}
$$

$$
(4.11) \qquad \widehat{\theta}_j - \theta_j = \int \Delta\psi_j\psi_j + o_p(n^{-1/2}).
$$

*Step* (iii): *Approximations to* $\Delta$. We may Taylor-expand $X_i(T_{ij})$ about $X_i(u)$, obtaining

$$
(4.12) \qquad X_i(T_{ij}) = X_i(u) - U_{ij}X_i'(u) + \tfrac{1}{2}U_{ij}^2 X''(u_{ij}),
$$

where $U_{ij} = u - T_{ij}$ and the random variable $u_{ij}$ lies between $u$ and $T_{ij}$ and is of course independent of the errors $\varepsilon_{rs}$. For given $u, v \in \mathscr{I}$, define $Z^{[1]}_{ijk} = \{X_i(u) + \varepsilon_{ij}\}\{X_i(v) + \varepsilon_{ik}\}$, $V_{ik} = v - T_{ik}$, $Z^{[2]}_{ijk} = U_{ij}X_i'(u)X_i(v) + V_{ik}X_i(u)X_i'(v)$ and $Z^{[3]}_{ijk} = U_{ij}X_i'(u)\varepsilon_{jk} + V_{ik}X_i'(u)\varepsilon_{jk}$. Let $\widehat{\phi}^{[\ell]}$ denote the version of $\widehat{\phi}$ obtained on replacing $Z_{ijk}$ by $Z^{[\ell]}_{ijk}$. In the following calculations, we may set $\mu \equiv 0$, without loss of generality. Using (4.12), and its analogue for expansion about $X_i(v)$ rather than $X_i(u)$, we may write

$$
\begin{aligned}
(4.13) \quad Y_{ij}Y_{ik} &= \{X_i(u) - U_{ij}X_i'(u) + \tfrac{1}{2}U_{ij}^2 X_i''(u_{ij}) + \varepsilon_{ij}\} \\
&\quad \times \{X_i(v) - V_{ik}X_i'(v_{ij}) + \tfrac{1}{2}V_{ij}^2 X_i''(v) + \varepsilon_{ik}\} \\
&= Z_{ijk}^{[1]} - Z_{ijk}^{[2]} - Z_{ijk}^{[3]} + Z_{ijk}^{[4]},
\end{aligned}
$$

where $Z_{ijk}^{[4]}$ is defined by (4.13). Using standard arguments for deriving uniform convergence rates it may be proved that for some $\eta > 0$, and under the bandwidth conditions for either part of Theorem 1,

$$
\sup_{(u,v)\in\mathscr{I}^2} |\widehat{\phi}^{[4]}(u,v) - E'\{\widehat{\phi}^{[4]}(u,v)\}| = O_p(n^{-(1/2)-\eta}).
$$

[Note that the data $Z_{ijk}^{[4]}$, from which $\widehat{\phi}^{[4]}$ is computed, contain only quadratic terms in $(U_{ij}, V_{ik})$. When the kernel weights are applied for constructing $\widehat{\phi}^{[4]}$, only triples $(i,j,k)$ for which $|U_{ij}|, |V_{ij}| \leq \mathrm{const.}\, h_\phi$ make a nonvanishing contribution to the estimator. This fact ensures the relatively fast rate of convergence.] Therefore, uniformly on $\mathscr{I}^2$, and under either set of bandwidth conditions,

$$
\begin{aligned}
(4.14) \quad \widehat{\phi} - E'\widehat{\phi} &= \widehat{\phi}^{[1]} - E'\widehat{\phi}^{[1]} - (\widehat{\phi}^{[2]} - E'\widehat{\phi}^{[2]}) - (\widehat{\phi}^{[3]} - E'\widehat{\phi}^{[3]}) \\
&\quad + o_p(n^{-1/2}).
\end{aligned}
$$

Put $Y_{ij}^{[1]}(u) = X_i(u) + \varepsilon_{ij}$, and let $\widehat{\mu}^{[1]}$ denote the version of $\widehat{\mu}$ obtained on replacing $Y_{ij}$ by $Y_{ij}^{[1]}$. Define $\widehat{\mu}^{[2]}$ by $\widehat{\mu} = \widehat{\mu}^{[1]} + \widehat{\mu}^{[2]}$. Conventional arguments for deriving uniform convergence rates may be employed to show that for some $\eta > 0$, and under either set of bandwidth conditions,

$$
(4.15) \quad \sup_{u\in\mathscr{I}} |\widehat{\mu}^{[1]}(u) - E'\{\widehat{\mu}^{[1]}(u)\}| = O_p(n^{-(1/4)-\eta}),
$$

$$
(4.16) \quad \sup_{u\in\mathscr{I}} |\widehat{\mu}^{[2]}(u) - E'\{\widehat{\mu}^{[2]}(u)\}| = O_p(n^{-(1/4)-\eta}),
$$

where (4.15) makes use of the property that $h_\mu \geq n^{\eta-(1/2)}$ for some $\eta > 0$. Combining (4.14)–(4.16) we deduce that, for either set of bandwidth conditions,

$$
(4.17) \quad \Delta = \Delta_0 + \Delta_1 - \Delta_2 - \Delta_3 - \Delta_4 - \Delta_5 + \Delta_6,
$$

where $\Delta_k = \widehat{\phi}^{[k]} - E'\widehat{\phi}^{[k]}$ for $k = 1, 2, 3$,

$$
\begin{aligned}
\Delta_0(u,v) &= E'\{\widehat{\phi}(u,v)\} - \{E'\widehat{\mu}(u)\}\{E'\widehat{\mu}(v)\} - \psi(u,v), \\
\Delta_4(u,v) &= E'\{\widehat{\mu}(u)\}\{\widehat{\mu}^{[1]}(v) - E'\widehat{\mu}^{[1]}(v)\},
\end{aligned}
$$

$\Delta_5(u, v) = \Delta_4(v, u)$ and $\sup_{(u,v) \in \mathscr{I}^2} |\Delta_6(u, v)| = o_p(n^{-1/2})$. Direct calculation may be used to show that for each $r_0, s_0 \geq 1$, and for either set of bandwidth conditions,

$$\max_{1 \leq k \leq 3} \max_{1 \leq r, s \leq r_0} E' \left\{ \iint_{\mathscr{I}^2} \Delta_k(u, v) \psi_r(u) \psi_s(v) \, du \, dv \right\}^2 = O_p(n^{-1}),$$

$$\max_{k=2,3} \max_{1 \leq r, s \leq r_0} E' \left\{ \iint_{\mathscr{I}^2} \Delta_k(u, v) \psi_r(u) \psi_s(v) \, du \, dv \right\}^2 = o_p(n^{-1}),$$

$$\max_{1 \leq k \leq 3} \sup_{s \geq 1} E' \left\{ \iint_{\mathscr{I}^2} \Delta_k(u, v) \psi_r(u) \psi_s(v) \, du \, dv \right\}^2 = O_p(n^{-1}),$$

$$E'(\|\Delta_1\|_{(j)}^2) = O_p\{(nh_\phi)^{-1}\},$$

$$\max_{k=2,3} E'(\|\Delta_k\|_{(j)}^2) = O_p(n^{-1}),$$

$$\max_{k=1,2} \max_{r \geq 1, 1 \leq s \leq s_0} E' \left[ \iiint_{\mathscr{I}^2} E'\{\widehat{\mu}(u)\}\{\widehat{\mu}^{[k]}(v) - E'\widehat{\mu}^{[k]}(v)\} \right.$$
$$\left. \times \psi_r(u) \psi_s(v) \, du \, dv \right]^2 = O_p(n^{-1}),$$

$$\max_{k=1,2} \sup_{1 \leq r \leq r_0, s \geq 1} E' \left[ \iiint_{\mathscr{I}^2} E'\{\widehat{\mu}(u)\}\{\widehat{\mu}^{[k]}(v) - E'\widehat{\mu}^{[k]}(v)\} \right.$$
$$\left. \times \psi_r(u) \psi_s(v) \, du \, dv \right]^2 = O_p\{(nh_\mu)^{-1}\}.$$

Standard arguments show that

$$\sup_{(u,v) \in \mathscr{I}^2} |\Delta_0(u, v) - \psi(u, v)| = O_p(h^2).$$

Combining results from (4.17) down, and defining

$$(4.18) \qquad D_{j3} = \sum_{k: k \neq j} \{(\theta_j - \theta_k)^{-2} - \theta_j^{-2}\} \left( \int \Delta_0 \psi_j \psi_j \right)^2$$

[cf. (4.7)], we deduce from (4.10) and (4.11) that, for the bandwidth conditions in (a) and (b), respectively,

$$(4.19) \qquad \begin{aligned} \|\widehat{\psi}_j - \psi_j\|^2 &= D_{j3} + \theta_j^{-2} \|\Delta_0 + \Delta_1\|_{(j)}^2 - \theta_j^{-2} \left( \int \Delta_0 \psi_j \psi_j \right)^2 \\ &\quad + o_p\{(nh_\phi)^{-1} + h_\phi^4\} + O_p\{(nh_\mu)^{-3/2}\}, \end{aligned}$$

$$(4.20) \qquad \widehat{\theta}_j - \theta_j = \int (\Delta_0 + \Delta_1 - \Delta_4 - \Delta_5) \psi_j \psi_j + o_p(n^{-1/2}).$$

*Step* (iv): *Elucidation of bias contributions.* Let the function $\chi$ be as defined at (3.1). It may be proved that $E'\{\hat{\phi}(u,v)\} = \phi(u,v) + h_\phi^2\chi(u,v) + o_p(h_\phi^2)$, uniformly in $u,v \in \mathscr{I}$ with $|u-v| > \delta$ for any $\delta > 0$, and that $E'\{\hat{\phi}(u,v)\} = O_p(h_\phi^2)$, uniformly in $u,v \in \mathscr{I}$. Here one uses the fact that $E\{X(s)X(t)\} = \psi(s,t) + x(s)x(t)$; subsequent calculations involve replacing $(s,t)$ by $(T_{ij}, T_{ik})$ on the right-hand side.

Furthermore, $E'\{\hat{\mu}(u)\} = \mu(u) + O_p(h_\mu^2)$, uniformly in $u \in \mathscr{I}$. In view of these properties and results given in the previous paragraph, and noting that we assume $h_\mu = o(h_\phi)$ in the context of (4.19), and $h_\phi = o(n^{-1/4})$ for (4.20), we may replace $\Delta_0$ by $h_\phi^2\chi$ in the definition of $D_{j3}$ at (4.18), and in (4.19) and (4.20), without affecting the correctness of (4.19) and (4.20). Let $D_{j4}$ denote the version of $D_{j3}$ where $\Delta_0$ is replaced by $h_\phi^2\chi$.

Moment methods may be used to prove that

$$
\begin{aligned}
(4.21) \quad \|h_\phi^2\chi + \Delta_1\|_{(j)}^2 &= E'\|h_\phi^2\chi + \Delta_1\|_{(j)}^2 + o_p\{(nh_\phi)^{-1} + h_\phi^4\} \\
&= h_\phi^4\|\chi\|_{(j)}^2 + E'\|\Delta_1\|_{(j)}^2 + o_p\{(nh_\phi)^{-1} + h_\phi^4\}.
\end{aligned}
$$

Furthermore,

$$
\begin{aligned}
(4.22) \quad &D_{j4} + \theta_j^{-2}h_\phi^4\|\chi\|_{(j)}^2 - \theta_j^{-2}h_\phi^4\left(\int \chi\psi_j\psi_j\right)^2 \\
&= h_\phi^4 \sum_{k:\,k\neq j} (\theta_j - \theta_k)^{-2}\left(\int \chi\psi_j\psi_k\right)^2;
\end{aligned}
$$

compare (4.6). Combining (4.21) and (4.22) with the results noted in the previous paragraph, we deduce that, under the bandwidth conditions assumed for parts (a) and (b), respectively, of Theorem 1,

$$
\begin{aligned}
(4.23) \quad \|\hat{\psi}_j - \psi_j\|^2 &= \theta_j^{-2}E'\|\Delta_1\|_{(j)}^2 + h_\phi^4 \sum_{k:k\neq j}(\theta_j - \theta_k)^{-2}\left(\int \chi\psi_j\psi_k\right)^2 \\
&\quad + o_p\{(nh_\phi)^{-1} + h_\phi^4\} + O_p\{(nh_\mu)^{-3/2}\},
\end{aligned}
$$

$$
(4.24) \quad \hat{\theta}_j - \theta_j = \int (\Delta_1 - 2\Delta_4)\psi_j\psi_j + o_p(n^{-1/2}).
$$

*Step* (v): *Calculation of $E'\|\Delta_1\|_{(j)}^2$.* Since $E'\|\Delta_1\|_{(j)}^2 = \int_{\mathscr{I}} \xi_1(u)\,du$, where

$$
\xi_1(u) = \iint_{\mathscr{I}^2} \xi_2(u,v_1,v_2)\psi_j(v_1)\psi_j(v_2)\,dv_1\,dv_2
$$

and $\xi_2(u,v_1,v_2) = E'\{\Delta_1(u,v_1)\Delta_1(u,v_2)\}$, then we shall compute $\xi_2(u,v_1,v_2)$.

Recall the definition of $\hat{\phi}$ at (2.5). An expression for $\Delta_1$ is the same, except that we replace $R_{rs}$ in the formula for $\hat{\phi}$ by $Q_{rs}$, say, which is defined by replacing $Z_{ijk}$, in the definition of $R_{rs}$ in Section 2, by $Z_{ijk}^{[1]} - E'(Z_{ijk}^{[1]})$. It can then be seen that if, in the formulae in the previous paragraph, we replace $\xi_2(u, v_1, v_2)$ by $\xi_3(u, v_1, v_2) = E'\{\Delta^*(u, v_1)\Delta^*(u, v_2)\}$, where $\Delta^*(u, v) = A_1 Q_{00}/B$ and $A_1$ and $B$ are as in Section 2, then we commit an error of smaller order than $(nh_\phi)^{-1} + h_\phi^4$ in the expression for $E'\|\Delta_1\|_{(j)}^2$.

Define $x = X - \mu$ and $\beta(s_1, t_1, s_2, t_2) = E\{x(s_1)x(t_1)x(s_2)x(t_2)\}$. In this notation,

$$
\begin{aligned}
&B(u, v_1)B(u, v_1)\xi_3(u, v_1, v_2)/A_1(u, v_1)A_1(u, v_2) \\
&= \sum_{i=1}^n \sum_{j_1<k_1} \sum_{j_2<k_2} \beta(T_{ij_1}, T_{ik_1}, T_{ij_2}, T_{ik_2}) \\
&\hspace{3em} \times K\left(\frac{T_{ij_1} - u}{h_\phi}\right)K\left(\frac{T_{ik_1} - v_1}{h_\phi}\right) \\
&\hspace{3em} \times K\left(\frac{T_{ij_2} - u}{h_\phi}\right)K\left(\frac{T_{ik_2} - v_2}{h_\phi}\right) \\
&\quad + \sigma^2 \sum_{i=1}^n \sum_{j<k} K\left(\frac{T_{ij} - u}{h_\phi}\right)^2 K\left(\frac{T_{ik} - v_1}{h_\phi}\right)K\left(\frac{T_{ik} - v_2}{h_\phi}\right).
\end{aligned}
\tag{4.25}
$$

Contributions to the fivefold series above, other than those for which $(j_1, k_1) = (j_2, j_2)$, make an asymptotically negligible contribution. Omitting such terms, the right-hand side of (4.25) becomes

$$
\sum_{i=1}^n \sum_{j<k} \{\beta(T_{ij}, T_{ik}, T_{ij}, T_{ik}) + \sigma^2\} K\left(\frac{T_{ij} - u}{h_\phi}\right)^2 K\left(\frac{T_{ik} - v_1}{h_\phi}\right)K\left(\frac{T_{ik} - v_2}{h_\phi}\right).
$$

Multiplying by $\psi_j(v_1)\psi_j(v_2)A_1(u, v_1)A_1(u, v_2)\{B(u, v_1)B(u, v_2)\}^{-1}$, integrating over $u, v_1, v_2$, and recalling that $N = \sum_i m_i(m_i - 1)$, we deduce that

$$
\begin{aligned}
E'\|\Delta_1\|_{(j)}^2 &\sim (Nh_\phi^4)^{-1}\int_{\mathscr{I}} f(u)^{-2}\, du \\
&\quad \times \iiiint_{\mathscr{I}^4} \{\beta(t_1, t_2, t_1, t_2) + \sigma^2\} \\
&\hspace{3em} \times K\left(\frac{t_1 - u}{h_\phi}\right)^2 K\left(\frac{t_2 - v_1}{h_\phi}\right)K\left(\frac{t_2 - v_2}{h_\phi}\right)\psi_j(v_1)\psi_j(v_2) \\
&\hspace{3em} \times \{f(v_1)f(v_2)\}^{-1} f(t_1)f(t_2)\, dt_1\, dt_2\, dv_1\, dv_2 \\
&\sim (Nh_\phi)^{-1}\int \cdots \int \{f(t_1)f(t_2)\}^{-1}\{\beta(t_1, t_2, t_1, t_2) + \sigma^2\} \\
&\hspace{3em} \times K(s)^2 K(s_1)K(s_2)\psi_j(t_2)^2\, dt_1\, dt_2\, ds_1\, ds_2\, ds \\
&= (Nh_\phi)^{-1}C_1.
\end{aligned}
$$

Result (3.3) follows from this property and (4.23).

*Step* (vi): *Limit distribution of* $Z_j \equiv \int (\Delta_1 - 2\Delta_4)\psi_j\psi_j$. It is straightforward to prove that the vector $Z$, of which the $j$th component is $Z_j$, is asymptotically normally distributed with zero mean. We conclude our proof of part (b) of Theorem 1 by finding its asymptotic covariance matrix, which is the same as the limit of the covariance conditional on the set of observation times $T_{ij}$. In this calculation we may, without loss of generality, take $\mu \equiv 0$.

Observe that

$$(4.26) \quad \operatorname{cov}'(Z_r, Z_s) = c_{11}(r,s) - 2c_{14}(r,s) - 2c_{14}(s,r) + 4c_{44}(r,s),$$

where the dash in $\operatorname{cov}'$ denotes conditioning on observation times, and

$$c_{ab}(r,s) = \iiiint_{\mathscr{I}^4} E'\{\Delta_a(u_1, v_1)\Delta_b(u_2, v_2)\}$$
$$\times \psi_r(u_1)\psi_r(v_1)\psi_s(u_2)\psi_s(v_2)\, du_1\, dv_1\, du_2\, dv_2.$$

We shall compute asymptotic formulae for $c_{11}$, $c_{14}$ and $c_{44}$.

Pursuing the argument in step (iv) we may show that $E'(\Delta_1\Delta_1)$ is asymptotic to

$$\frac{A_1(u_1, v_1)A_1(u_2, v_2)}{B(u_1, v_1)B(u_2, v_1)}$$
$$\times \left\{\sum_{i=1}^{n}\sum_{j_1 < k_1}\sum_{j_2 < k_2}\beta(T_{ij_1}, T_{ik_1}, T_{ij_2}, T_{ik_2})\right.$$
$$\times K\left(\frac{T_{ij_1} - u_1}{h_\phi}\right)K\left(\frac{T_{ik_1} - v_1}{h_\phi}\right)$$
$$\times K\left(\frac{T_{ij_2} - u_2}{h_\phi}\right)K\left(\frac{T_{ik_2} - v_2}{h_\phi}\right)$$
$$+ \sigma^2\sum_{i=1}^{n}\sum_{j<k}\sum K\left(\frac{T_{ij} - u_1}{h_\phi}\right)K\left(\frac{T_{ij} - u_2}{h_\phi}\right)$$
$$\left.\times K\left(\frac{T_{ik} - v_1}{h_\phi}\right)K\left(\frac{T_{ik} - v_2}{h_\phi}\right)\right\}.$$

The ratio $A_1A_1/BB$ to the left is asymptotic to $\{S_{00}(u_1, v_1)S_{00}(u_2, v_2)\}^{-1}$, and so to $\{N^2 h_\phi^4 f(u_1)f(u_2)f(v_1)f(v_2)\}^{-1}$. Therefore,

$$c_{11}(r,s) \overset{p}{\sim} N^{-2}\left[\sum_{i=1}^{n}\sum_{j_1 < k_1}\sum_{j_2 < k_2}\{f(T_{ij_1})f(T_{ik_1})f(T_{ij_2})f(T_{ik_2})\}^{-1}\right.$$

$$\times \beta(T_{ij_1}, T_{ik_1}, T_{ij_2}, T_{ik_2})$$

$$\times \psi_r(T_{ij_1}) \psi_r(T_{ik_1}) \psi_s(T_{ij_2}) \psi_s(T_{ik_2})$$

$$(4.27) \qquad \times \iiiint_{\mathscr{I}^4} K(w_1) K(w_2) K(w_3) K(w_4)\, dw_1 \cdots dw_4$$

$$+ \sigma^2 \sum_{i=1}^n \sum_{j<k} \sum \{f(T_{ij}) f(T_{ik})\}^{-2} \psi_r(T_{ij}) \psi_r(T_{ik}) \psi_s(T_{ij}) \psi_s(T_{ik})$$

$$\times \iiiint_{\mathscr{I}^4} K(w_1) K(w_2) K(w_3) K(w_4)\, dw_1 \cdots dw_4 \Bigg]$$

$$\overset{p}{\sim} N^{-2} \{\nu(r,s) + N\sigma^2 c(r,s)^2\}.$$

Observe next that, defining $\gamma(u,v,w) = E\{X(u)X(v)X(w)\} - \phi(u,v)\mu(w)$, it may be shown that $S_0(v_2) S_{00}(u_1,v_1) E'\{\Delta_1(u_1,v_1)\Delta_4(u_2,v_2)\}/\mu(u_2)$ is asymptotic to

$$\sum_{i=1}^n \sum_{j_1<k_1} \sum \sum_{j_2=1}^{m_i} W_{ij_1k_1}(u_1,v_1) W_{ij_2}(v_2)$$

$$\times E'([\{X_i(u_1) + \varepsilon_{ij_1}\}\{X_i(v_1) + \varepsilon_{ik_1}\} - \phi(u_1,v_1)]$$

$$\times \{X_i(v_2) + \varepsilon_{ij_2} - \mu(v_2)\})$$

$$= \sum_{i=1}^n \sum_{j_1<k_1} \sum \sum_{j_2=1}^{m_i} W_{ij_1k_1}(u_1,v_1) W_{ij_2}(v_2)$$

$$\times \{\gamma(u_1,v_1,v_2) + \sigma^2 \delta_{j_1j_2}\mu(v_1) + \sigma^2 \delta_{j_2k_1}\mu(u_1)\}.$$

The terms in $\sigma^2$ can be shown to make asymptotically negligible contributions to $c_{14}$; the first of them vanishes unless $u_1$ is within $O(h_\phi)$ of $v_2$, and the second unless $v_1$ is within $O(h_\phi)$ of $v_2$. Furthermore, defining $N_1 = N_1(n) = \sum_{i\leq n} m_i$, we have $S_0(v_2) \sim_p N_1 h_\mu f(v_2)$ and $S_{00}(u_1,v_1) \sim_p N h_\phi^2 f(u_1) f(v_1)$. From these properties, and the fact that we are assuming (without loss of generality) that $\mu \equiv 0$, it may be proved that

$$(4.28) \qquad |c_{14}(r,s)| + |c_{44}(r,s)| = o_p(n^{-1}).$$

Results (4.24) and (4.26)–(4.28) imply that the covariance matrix in the central limit theorem for the vector of values of $\hat{\theta}_j - \theta_j$ has the form stated in part (b) of Theorem 1.

**Acknowledgments.** We wish to thank the Editor, an Associate Editor and two reviewers for helpful comments that led to several improvements in the paper.

P. HALL
CENTER FOR MATHEMATICS AND ITS APPLICATIONS
AUSTRALIAN NATIONAL UNIVERSITY
CANBERRA ACT 0200
AUSTRALIA

H.-G. MÜLLER
J.-L. WANG
DEPARTMENT OF STATISTICS
UNIVERSITY OF CALIFORNIA
ONE SHIELDS AVENUE
DAVIS, CALIFORNIA 95616
USA
E-MAIL: mueller@wald.ucdavis.edu